\documentclass[10pt]{amsart}
\usepackage{amssymb,mathrsfs}
\usepackage{amsmath}
\usepackage{enumitem}
\usepackage{tikz-cd}
\usepackage[OT2,T1]{fontenc}
\DeclareSymbolFont{cyrletters}{OT2}{wncyr}{m}{n}
\DeclareMathSymbol{\Sha}{\mathalpha}{cyrletters}{"58}
\newcommand{\smallpagebreak}{{\par\vspace{2 mm}\noindent}}
\newcommand{\medpagebreak}{{\par\vspace{4 mm}\noindent}}

\renewcommand\a{\alpha}

\newcommand\G{\Gamma}

\newcommand\la{\lambda}

\newcommand\lrt{\longrightarrow}

\newcommand\psb{\par\smallpagebreak}

\newcommand\BZ{\mathbb{Z}}
\newcommand\BC{\mathbb{C}}
\newcommand\BR{\mathbb{R}}
\newcommand\BQ{\mathbb{Q}}
\newcommand\BH{\mathbb{H}}

\newcommand\fg{\mathfrak{g}}

\newcommand\pmm{\par\medpagebreak}
\theoremstyle{definition}
\newtheorem{theorem}{Theorem}[section]
\newtheorem{definition}[theorem]{Definition}
\newtheorem{corollary}[theorem]{Corollary}
\newtheorem{example}{Example}[section]
\newtheorem{result}{Result}
\newtheorem{conjecture}{Conjecture}

\newtheorem{remark}[theorem]{Remark}

\newtheorem{proposition}[theorem]{Proposition}
\numberwithin{equation}{section}

\begin{document}
\title[Remarks on The Birch-Swinnerton-Dyer Conjecture]{\large Remarks on The Birch-Swinnerton-Dyer Conjecture}
\begin{abstract}
We give a brief description of the Birch-Swinnerton-Dyer
conjecture which is one of the seven Clay problems and present several related conjectures. We describe the relation between the nilpotent orbits of $SL(2,\BR)$ and CM points.
\end{abstract}

\author{Jae-Hyun Yang}
\address{Department of Mathematics, Inha University,Incheon 22212, Korea}
\email{jhyang@inha.ac.kr }
\maketitle
\vskip 0.5cm
\section{Introduction}    
\vskip 0.5cm \ \ \ On May 24, 2000, the Clay Mathematics Institute
(CMI for short) announced that it would award prizes of 1 million
dollars each for solutions to seven mathematics problems. These
seven problems are

\begin{enumerate}[label=Problem \arabic{enumi}.,itemsep=6pt,leftmargin=*]
\item The ``P versus NP" Problem
\item The Riemann Hypothesis
\item The Poincar{\'e} Conjecture
\item The Hodge Conjecture
\item The Birch-Swinnerton-Dyer Conjecture (briefly, the BSD conjecture)
\item The Navier-Stokes Equations\,: Prove or disprove the existence and \\
		smoothness of solutions to the three dimensional Navier-Stokes equations.
\item Yang-Mills Theory : Prove that quantum Yang-Mills fields exist  \\
		and have a mass gap.
\end{enumerate}
\par Problem 1 is arisen from theoretical computer science, Problem 2 and Problem 5
from number theory, Problem 3 from topology, Problem 4 from
algebraic geometry and topology, and finally problem 6 and 7 are
related to physics. For more details on some stories about these
problems, we refer to Notices of AMS, vol. 47, no. 8,
pp.\,877-879\,(September 2000) and the homepage of CMI. In 2003,
Problem 3 was solved by Grisha Perelman \cite{P1, P2, P3}. We refer to \cite{CZ1, CZ2, KL, MT1, MT2} for more details on Perelman's work. Recently Bhargava's school computed the Selmer groups of an elliptic curve and so solved Problem 5 partially.

\vskip 2mm
The purpose of this paper is to describe the relation between the nilpotent orbits of $SL(2,\BR)$ and CM points and to present several conjectures relating to the BSD Conjecture.

\vskip 2mm
The paper is organized as follows.
From Section 2 to Section 5, we will explain Problem 5, that is, the
Birch-Swinnerton-Dyer conjecture which was proposed by the English mathematicians, B. Birch and H. P. F. Swinnerton-Dyer \cite{BSD1, BSD2} around 1960s
in some detail. This conjecture
says that if $E$ is an elliptic curve defined over $\BQ$, then the
algebraic rank of $E$ equals the analytic rank of $E.$ In 2001,
the Shimura-Taniyama conjecture stating that any elliptic curve defined over $\BQ$
is modular was shown to be true by Breuil,
Conrad, Diamond and Taylor\,\cite{BCDT}. This fact shed some lights on the solution of the BSD conjecture.
In Section 6, we describe the connection between the
heights of Heegner points on modular curves $X_0(N)$ and Fourier
coefficients of modular forms of half integral weight or of the
Jacobi forms corresponding to them by the Skoruppa-Zagier
correspondence. Most of the materials in Section 2--6 were already printed in \cite{Y5}. In Section 7, we briefly review the works done recently by the school of Manjul Bhargava \cite{BS1, BS2, BS3, BS4, BSk, BES}. In Section 8, we describe the adjoint orbits of $SL(2,\BR)$ in its Lie algebra $\frak{s}\frak{l} (2,\BR)$ explicitly. In the final section, we describe the relation between the nilpotent orbits of $SL(2,\BR)$ and CM points. We
propose several conjectures relating to the Birch-Swinnerton-Dyer conjecture.

\pmm\noindent
{\bf Notations\,:} We denote by $\BQ,\, \BR$ and $\BC$ the fields
of rational numbers, real numbers and complex numbers
respectively. $\BZ$ and $\BZ^+$ denotes the ring of integers and
the set of positive integers respectively. $\mathbb H$ denotes the Poincar{\'e} upper half plane.

\vskip 8mm
\section{The Mordell-Weil Group}
\vskip 0.3cm
\ \ \ A curve $E$ is said to be an {\it elliptic
curve} over $\BQ$ if it is a nonsingular projective curve of genus
1 with its affine model
\begin{equation}
y^2 =f(x),
\end{equation}
 where $f(x)$ is a
polynomial of degree 3 with integer coefficients and with 3
distinct roots over $\BC$. An elliptic curve over $\BQ$ has an
abelian group structure with distinguished element $\infty$ as an
identity element. The set $E(\BQ)$ of rational points given by
\begin{equation}
E(\BQ)=\left\{\,(x,y)\in\BQ^2\,|\ y^2=f(x)\,\right\}\cup
\{\infty\}
\end{equation}also has an abelian group
structure.\pmm L. J. Mordell\,(1888-1972) \cite{Mo} proved the
following theorem in 1922.\psb\noindent

\begin{theorem}[Mordell,1922]
 $E(\BQ)$ is finitely generated, that is,
 \begin{equation*}
  E(\BQ)\cong \BZ^r \oplus E_{\text{tor}}(\BQ),
 \end{equation*}
where $r$ is a nonnegative
integer and $E_{\text{tor}}(\BQ)$ is the torsion subgroup of
$E(\BQ)$.
\end{theorem}
\vskip 0.1cm
\begin{definition}
	Around 1930, A. Weil\,(1906-1998) proved that the set $A(\BQ)$ of rational points on an
	abelian variety $A$ defined over $\BQ$ is finitely generated. An
	elliptic curve is an abelian variety of dimension one. Therefore
	$E(\BQ)$ is called the {\it Mordell}-{\it Weil group} and the
	integer $r$ is said to be the {\sf algebraic rank} of $E$.
\end{definition}
\vskip 0.1cm
In 1977, B. Mazur\,(1937-\ )\cite{Ma1} discovered the structure of the torsion subgroup $E_{\text{tor}}(\BQ)$ completely using a deep
theory of modular curves.
\begin{theorem}[Mazur, 1977] Let $E$ be an elliptic curve defined over
$\BQ$. Then the torsion subgroup $E_{\text{tor}}(\BQ)$ is
isomorphic to the following 15 groups $$\BZ/n\BZ\quad (1\leq n\leq
10,\ n=12),$$ $$\BZ/2\BZ \times \BZ/2n\BZ\quad (1\leq n\leq 4).$$
\end{theorem}
\vskip 0.1cm
\indent E. Lutz\,(1914-?) and T. Nagell\,(1895-?) obtained the
following result independently.
\begin{theorem}[Lutz, 1937;\ Nagell, 1935]
Let $E$ be an elliptic curve
defined over $\BQ$ given by $$E\,:\quad y^2=x^2+ax+b,\quad
a,b\in\BZ,\ 4a^3+27b^2\not= 0.$$ Suppose that $P=(x_0,y_0)$ is an
element of the torsion subgroup $E_{\text{tor}}(\BQ)$. Then \psb
(a) $\quad x_0,y_0\in \BZ$, and \psb (b) $2P=0$\ \ or\ \ $y_0^2 |
(4a^3+27b^2).$
\end{theorem}
\vskip 0.1cm
We observe that the above theorem gives an
effective method for bounding $E_{\text{tor}}(\BQ)$. According to
Theorem B and C, we know the torsion part of $E(\BQ)$
satisfactorily. But we have no idea of the free part of $E(\BQ)$
so far. As for the algebraic rank $r$ of an elliptic curve $E$
over $\BQ$, Noam Elkies found an example of an elliptic curve of rank 28 in 2006. Indeed, that elliptic curve is given by
$$E_e:\ \ y^2+xy+y=x^3-x^2-\alpha \,x+\beta$$ has its
algebraic rank 28. Here
\begin{equation*}
\alpha= 20067762415575526585033208209338542750930230312178956502
\end{equation*}
and
\begin{equation*}
\beta=34481611795030556467032985690390720374855944359319180361266008296291939448732243429.
\end{equation*}
Elkies also computed 28 generators of $E_e(\BQ)$
(cf. https://web.math.pmf.unizg.hr/~duje/tors/rk28.html).

\pmm\noindent
\begin{conjecture}
 Given a nonnegative integer $n$, there is an elliptic curve $E$ over $\BQ$ with its algebraic rank $n$.
\end{conjecture}
\vskip 0.1cm
The algebraic rank of an
elliptic curve is an invariant under an isogeny. Here an isogeny
of an elliptic curve $E$ means a holomorphic map $\varphi:
E(\BC)\lrt E(\BC)$ satisfying the condition $\varphi(0)=0.$

\vskip 8mm
\section{Modular Elliptic Curves}
\vskip 0.25cm \ \ \ For a positive integer $N\in\BZ^+,$ we let
$$\Gamma_0(N):=\left\{ \begin{pmatrix}a & b\\ c & d\end{pmatrix} \in
SL(2,\BZ)\,\big| \ N|c\ \right\}$$
be the Hecke subgroup of
$SL(2,\BZ)$ of level $N$. Let $\BH$ be the Poincar{\'e} upper half plane. Then
$$Y_0(N)=\Gamma_0(N)\backslash \BH$$
is a noncompact surface, and
\begin{equation}
X_0(N)=\Gamma_0(N)\backslash\BH\cup \BQ\cup \{\infty\}
\end{equation} is a
compactification of $Y_0(N).$ We recall that a {\it cusp form} of
weight $k\geq 1$ and level $N\geq 1$ is a holomorphic function $f$
on $\BH$ such that for all $\begin{pmatrix}a & b\\ c & d\end{pmatrix}  \in
\Gamma_0(N)$ and for all $z\in\BH$, we have
$$f((az+b)/(cz+d))=(cz+d)^k f(z)$$ and $|f(z)|^2 (\text{Im}\,z)^k$
is bounded on $\BH$. We denote the space of all cusp forms of
weight $k$ and level $N$ by $S_k(N)$. If $f\in S_k(N)$, then it
has a Fourier expansion $$f(z)=\sum_{n=1}^{\infty}c_n(f) q^n,\quad
q:=e^{2\pi i z}$$ convergent for all $z\in\BH.$ We note that there
is no constant term due to the boundedness condition on $f$. Now
we define the $L$-series $L(f,s)$ of $f$ to be
\begin{equation}
L(f,s)=\sum_{n=1}^{\infty}c_n(f)\,n^{-s}.
\end{equation} \indent For
each prime $p\nmid N$, there is a linear operator $T_p$ on
$S_k(N)$, called the Hecke operator, defined by
$$(f|{T_p})(z)=p^{-1}\sum_{i=0}^{p-1}f((z+i)/p)+p^{k-1}(cpz+d)^k\cdot
f((apz+d)/(cpz+d))$$ for any $\begin{pmatrix}a & b\\ c & d\end{pmatrix}
\in SL(2,\BZ)$ with $c\equiv 0\, (N)$ and $d\equiv p\,(N).$ The
Hecke operators $T_p$ for $pN$ can be diagonalized on the
space $S_k(N)$ and a simultaneous eigenvector is called an {\it
eigenform.} If $f\in S_k(N)$ is an eigenform, then the
corresponding eigenvalues, $a_p(f)$, are algebraic integers and we
have $c_p(f)=a_p(f)\,c_1(f).$\psb

Let $\la$ be a place of the
algebraic closure ${\bar {\BQ}}$ in $\BC$ above a rational prime
$\ell$ and ${\bar {\BQ}}_{\la}$ denote the algebraic closure of
$\BQ_{\ell}$ considered as a ${\bar {\BQ}}$-algebra via $\la$. It
is known that if $f\in S_k(N),$ there is a unique continuous
irreducible representation
\begin{equation}
\rho_{f,\la}:\text{Gal}({\bar
{\BQ}}/\BQ)\lrt GL_2({\bar{\BQ}}_{\la})
\end{equation} such that for any
prime ${p\nmid N}  {\ell}$,\ $\rho_{f,\la}$ is unramified at $p$
and $\text{tr}\,\rho_{f,\la}(\text{Frob}_p)=a_p(f).$ The existence
of $\rho_{f,\la}$ is due to G. Shimura\,(1930- ) if $k=2$ \cite{Sh}, to
P. Deligne\,(1944- ) if $k>2$
\cite{D}
and to P. Deligne and J.-P.
Serre\,(1926- ) if $k=1$ \cite{DS}. Its irreducibility is due to Ribet
if $k>1$ \cite{R}, and to Deligne and Serre if $k=1$ \cite{DS}. Moreover
$\rho_{f,\la}$ is odd and potentially semi-stable at $\ell$ in the
sense of Fontaine. We may choose a conjugate of $\rho_{f,\la}$
which is valued in $GL_2(\mathcal{O}_{{\bar{\BQ}}_{\la}})$, and
reducing modulo the maximal ideal and semi-simplifying yields a
continuous representation \begin{equation}
{\bar\rho}_{f,\la}:\text{Gal}({\bar
{\BQ}}/\BQ)\lrt GL_2({\bar{\mathbb {F}}}_{\ell}),
\end{equation} which, up to
isomorphism, does not depend on the choice of conjugate of
$\rho_{f,\la}$. \pmm\noindent
\begin{definition}
	Let
	$\rho:\text{Gal}({\bar {\BQ}}/\BQ)\lrt GL_2({\bar{\BQ}}_{\ell})$
	be a continuous representation which is unramified outside
	finitely many primes and for which the restriction of $\rho$ to a
	decomposition group at $\ell$ is potentially semi-stable in the
	sense of Fontaine. We call $\rho$ {\it modular} if $\rho$ is
	isomorphic to $\rho_{f,\la}$ for some eigenform $f$ and some $\la	| \ell.$
\end{definition}
\vskip 1mm
\begin{definition} An elliptic curve $E$
defined over $\BQ$ is said to be {\it modular} if there exists a
surjective holomorphic map $\varphi: X_0(N)\lrt E(\BC)$ for some
positive integer $N$.
\end{definition}
\psb
In 2001 C. Breuil, B. Conrad, F.
Diamond and R. Taylor \cite{BCDT} proved that the Taniyama-Shimura
conjecture is true. \psb\noindent
\begin{theorem}[\cite{BCDT}, 2001]
An elliptic curve defined over $\BQ$ is modular.
\end{theorem}
\vskip 1mm

Let $E$ be
an elliptic curve defined over $\BQ$. For a positive integer
$n\in\BZ^+$, we define the isogeny $[n]:E(\BC)\lrt E(\BC)$ by
\begin{equation}
[n]P:=nP=P+\cdots+P\ (n\ \text{times}),\quad P\in
E(\BC).
\end{equation}For a negative integer $n$, we define the isogeny
$[n]:E(\BC)\lrt E(\BC)$ by $[n]P:=-[-n]P,\ P\in E(\BC)$, where
$-[-n]P$ denotes the inverse of the element $[-n]P$. And
$[0]:E(\BC)\lrt E(\BC)$ denotes the zero map. For an integer
$n\in\BZ,\ [n]$ is called the multiplication-by-$n$ homomorphism.
The kernel $E[n]$ of the isogeny $[n]$ is isomorphic to
$\BZ/n\BZ\oplus \BZ/n\BZ.$ Let
$$\text{End}(E)=\left\{\varphi:E(\BC)\lrt E(\BC),\ \text{an\
isogeny}\,\right\}$$ be the endomorphism group of $E$. An elliptic
curve $E$ over $\BQ$ is said to have {\it complex multiplication}
(or CM for short) if $$\text{End}(E)\not\subseteq \BZ\cong \left\{
[n]\vert\ n\in\BZ\,\right\},$$ that is, there is a nontrivial
isogeny $\varphi:E(\BC)\lrt E(\BC)$ such that $\varphi\not= [n]$
for all integers $n\in\BZ.$ Such an elliptic curve is called a CM
{\it curve}. For most of elliptic curves $E$ over $\BQ$, we have
$\text{End}(E)\cong \BZ.$

\vskip 10mm

\section{The $L$-Series of an Elliptic Curve}
\vskip 0.25cm \ \ \ Let $E$ be an elliptic curve over $\BQ$. The
$L$-series $L(E,s)$ of $E$ is defined  as the product of the local
$L$-factors\,:
\begin{equation}
L(E,s)=\prod_{p| \Delta_E} (1-a_p p^{-s})^{-1}
\cdot \prod_{{p \nmid} \Delta_E} (1-a_p
p^{-s}+p^{1-2s})^{-1},
\end{equation}where $\Delta_E$ is the
discriminant of $E$, $p$ is a prime, and if $p \nmid \Delta_E$,
$$a_p:=p+1-|{\bar E}(\mathbb{F}_p)|,$$ and if $p| \Delta_E,$ we
set $a_p:=0,\,1,\,-1$ if the reduced curve ${\bar E}/\mathbb{F}_p$
has a cusp at $p$, a split node at $p$, and a nonsplit node at $p$
respectively. Then $L(E,s)$ converges absolutely for
$\text{Re}\,s>{\frac 32}$ from the classical result that $|a_p|< 2
\sqrt{p}$ for each prime $p$ due to H. Hasse\,(1898-1971) and is
given by an absolutely convergent Dirichlet series. We remark that
$x^2-a_p x+p$ is the characteristic polynomial of the Frobenius
map acting on ${\bar E}(\mathbb{F}_p)$ by $(x,y)\mapsto (x^p,y^p).$
\pmm\noindent
\begin{conjecture}
	Let $N(E)$ be the conductor of
	an elliptic curve $E$ over $\BQ$\,(\cite{S},\,p.\,361). We set
	$$\Lambda(E,s):=N(E)^{s/2}\,(2\pi)^{-s}\,\Gamma(s)\, L(E,s),\quad
	\text{Re}\,s>{\frac 32}.$$ Then $\Lambda(E,s)$ has an analytic
	continuation to the whole complex plane and satisfies the
	functional equation $$\Lambda(E,s)=\epsilon\, \Lambda(E,2-s),\quad
	\epsilon=\pm 1.$$ \indent The above conjecture is now true because the Shimura-Taniyama conjecture is true\,(cf. Theorem E). We have
	some knowledge about analytic properties of $E$ by investigating
	the $L$-series $L(E,s)$. The order of $L(E,s)$ at $s=1$ is called
	the {\sf analytic rank} of $E$.
\end{conjecture}
\vskip 0.1cm
Now we explain the connection
between the modularity of an elliptic curve $E$, the modularity of
the Galois representation and the $L$-series of $E$. For a prime
$\ell$, we let $\rho_{E,\ell}$\,(resp. ${\bar\rho}_{E,\ell})$
denote the representation of $\text{Gal}({\bar\BQ}/\BQ)$ on the
$\ell$-adic Tate module\,(resp.\ the $\ell$-torsion) of
$E({\bar\BQ}).$ Let $N(E)$ be the conductor of $E$. Then it is
known that the following conditions are equivalent\,:
\vskip 5mm
\begin{enumerate}[itemsep=6pt]
	\item  The
	$L$-function $L(E,s)$ of $E$ equals the $L$-function $L(f,s)$ for
	some eigenform $f$.
	\item The $L$-function $L(E,s)$
	of $E$ equals the $L$-function $L(f,s)$ for some eigenform $f$
	\\ of weight 2 and level $N(E)$.
	\item For some
	prime $\ell$, the representation $\rho_{E,\ell}$ is modular.
	\item For all primes $\ell$, the representation $\rho_{E,\ell}$ is
	modular.
	\item There is a non-constant holomorphic map $X_0(N)\lrt E(\mathbb{C})$ for some
	\\positive integer $N$.
	\item There is a non-constant morphism $X_0(N(E))\lrt E$ which is
	defined over $\BQ.$
	\item  $E$ admits a hyperbolic
	uniformization of arithmetic type\,(cf.\,\cite{Ma2} and \cite{Y1}).
\end{enumerate}

\vskip 10mm
\section{The Birch-Swinnerton-Dyer conjecture}
\vskip 0.25cm \ \ \ Now we state the BSD conjecture.\psb\noindent
{\bf The BSD Conjecture.} Let $E$ be an elliptic curve over $\BQ$.
Then the algebraic rank of $E$ equals the analytic rank of $E.$
\pmm

I will describe some historical backgrounds about the BSD
conjecture. Around 1960, Birch\,(1931- ) and
Swinnerton-Dyer\,(1927- ) formulated a conjecture which determines
the algebraic rank $r$ of an elliptic curve $E$ over $\BQ$. The
idea is that an elliptic curve with a large value of $r$ has a
large number of rational points and should therefore have a
relatively large number of solutions modulo a prime $p$ on the
average as $p$ varies. For a prime $p$, we let $N(p)$ be the
number of pairs of integers $x,y\,(\text{mod}\,p)$ satisfying
(2.1) as a congruence (mod $p$). Then the BSD conjecture in its
crudest form says that we should have an asymptotic formula
\begin{equation}
\prod_{p<x}\dfrac{N(p)+1}{p}\ \sim C \ \,(\text{log}\,p)^r \quad
\text{as}\ x\lrt\infty
\end{equation}for some constant $C>0.$ If the
$L$-series $L(E,s)$ has an analytic continuation to the whole
complex plane\,(this fact is conjectured; cf.\, Conjecture F),
then $L(E,s)$ has a Taylor expansion  $$L(E,s)=c_0
(s-1)^m+c_1(s-1)^{m+1}+\cdots$$ at $s=1$ for some non-negative
integer $m\geq 0$ and constant $c_0\not= 0.$ The BSD conjecture
says that the integer $m$, in other words, the analytic rank of
$E$, should equal the algebraic rank $r$ of $E$ and furthermore
the constant $c_0$ should be given by
\begin{equation}
c_0=\lim_{s\rightarrow
1}{{L(E,s)}\over {(s-1)^m}}=\a\cdot R\cdot
|E_{\text{tor}}(\BQ)|^{-1}\cdot \Omega\cdot S,
\end{equation} where
$\a>0$ is a certain constant, $R$ is the elliptic regulator of
$E,\ |E_{\text{tor}}(\BQ)|$ denotes the order of the torsion
subgroup $E_{\text{tor}}(\BQ)$ of $E(\BQ)$, $\Omega$ is a simple
rational multiple (depending on the bad primes) of the elliptic
integral $$\int_{\gamma}^{\infty}{{dx}\over {\sqrt{f(x)}}}\qquad
(\gamma=\text{the\ largest\ root\ of}\ f(x)=0)$$ and $S$ is an
integer square which is supposed to be the order of the
Tate-Shafarevich group $\text{III}(E)$ of $E$. \psb

The Tate-Shafarevich group $\text{III}(E)$ of $E$ is a very interesting
subject to be investigated in the future. Unfortunately
$\text{III}(E)$ is still not known to be finite. So far an
elliptic curve whose Tate-Shafarevich group is infinite has not
been discovered. So many mathematicians propose the following.
\begin{conjecture}
	The Tate-Shafarevich group
	$\text{III}(E)$ of $E$ is finite.
\end{conjecture}
\vskip 1mm
There are some evidences
supporting the BSD conjecture. I will list these evidences
chronologically.
\begin{result}[Coates-Wiles\,\cite{CW},\,1977]
	Let $E$ be a CM curve over $\BQ$.
	Suppose that the analytic rank of $E$ is zero. Then the algebraic
	rank of $E$ is zero.
\end{result}
\begin{result}[Rubin\,\cite{R},\,1981] Let $E$ be a CM curve over $\BQ$. Assume
that the analytic rank of $E$ is zero. Then the Tate-Shafarevich
group $\text{III}(E)$ of $E$ is finite.
\end{result}
\begin{result}[Gross-Zagier\,\cite{GZ},\,1986\,;\ \cite{BCDT},\,2001] Let $E$ be an
elliptic curve over $\BQ$. Assume that the analytic rank of $E$ is
equal to one and $\epsilon=-1$\,(cf.\,Conjecture F). Then the
algebraic rank of $E$ is equal to or bigger than one.
\end{result}
\begin{result}[Gross-Zagier\,\cite{GZ},\,1986] There
exists an elliptic curve $E$ over $\BQ$ such that
$\text{rank}\,E(\BQ)=\text{ord}_{s=1} L(E,s)=3$. For instance, the
elliptic curve ${\tilde E}$ given by $${\tilde E}\ :\quad
-139\,y^2=x^3+10\,x^2-20\,x+8$$ satisfies the above property.
\end{result}
\begin{result}[Kolyvagin\,\cite{K2},\,1990\,:\,Gross-Zagier\,\cite{GZ},\,1986\,:\,Bump-Friedberg-Hoffstein\,\cite{BFH},\,
	\\1990\,:\,Murty-Murty\,\cite{MM},\,1990\,:\,\cite{BCDT},\,2001] Let $E$ be
an elliptic curve over $\BQ$. Assume that the analytic rank of $E$
is 1 and $\epsilon=-1$. Then algebraic rank of $E$ is equal to 1.
\end{result}
\begin{result}[Kolyvagin\,\cite{K2},\,1990\,:\,Gross-Zagier\,\cite{GZ},\,1986\,:\,Bump-Friedberg-Hoffstein\,\cite{BFH}, \\
\,1990\,:\,Murty-Murty\,\cite{MM},\,1990\,:\,\cite{BCDT},\,2001] Let $E$ be
an elliptic curve over $\BQ$. Assume that the analytic rank of $E$
is zero and $\epsilon=1$. Then algebraic rank of $E$ is equal to
zero.
\end{result}
 Cassels proved the fact that if an elliptic curve over
$\BQ$ is isogeneous to another elliptic curve $E'$ over $\BQ$,
then the BSD conjecture holds for $E$ if and only if th e BSD
conjecture holds for $E'$.

\vskip 1cm
\section{Jacobi Forms and Heegner Points}
\vskip 0.25cm
\indent In this section, I shall describe the result of
Gross-Kohnen-Zagier\,\cite{GKZ} roughly.
\psb \indent
First we begin with giving
the definition of Jacobi forms. By definition a Jacobi form of
weight $k$ and index $m$ is a holomorphic complex valued function
$\phi(z,w)\,(z\in \BH,\,z\in\BC)$ satisfying the transformation
formula
\begin{align} \phi\left( {\frac{az+b}{cz+d}}, {{w+\la
z+\mu}\over{cz+d}}\right)=&e^{-2\pi i\left\{ cm(w+\la z+\mu)^2
(cz+d)^{-1}-m(\la^2 z+2\la w)\right\}} \\ &\ \times (cz+d)^k
\,\phi(z,w) \nonumber\end{align}
for all $\begin{pmatrix}a & b\\ c & d\end{pmatrix} \in SL(2,\BZ)$ and $(\la,\mu)\in\BZ^2$, and having a
Fourier expansion of the form \begin{equation}
\phi(z,w)=\sum_{\scriptstyle
n,r\in \BZ^2 \atop\scriptstyle r^2\leq 4mn} c(n,r)\,e^{2\pi
i(nz+rw)}.
\end{equation}
\indent We remark that the Fourier coefficients
$c(n,r)$ depend only on the discrimnant $D=r^2-4mn$ and the
residue $r\,(\text{mod}\ 2m).$ From now on, we put
$\Gamma_1:=SL(2,\BZ).$ We denote by $J_{k,m}(\Gamma_1)$ the space
of all Jacobi forms of weight $k$ and index $m$. It is known that
one has the following isomorphisms
\begin{equation}
[\Gamma_2,k]^M\cong
J_{k,1}(\G_1)\cong M_{k-{\frac 12}}^+(\G_0(4))\cong
[\G_1,2k-2],
\end{equation} where $\G_2$ denotes the Siegel modular
group of degree 2, $[\G_2,k]^M$ denotes the Maass space introduced
by H. Maass\,(1911-1993)\,(cf.\,\cite{M1, M2, M3}), $M_{k-{\frac
12}}^+(\G_0(4))$ denotes the Kohnen space introduced by W.
Kohnen\,\cite{Koh} and $[\G_1,2k-2]$ denotes the space of modular forms
of weight $2k-2$ with respect to $\G_1$. We refer to \cite{Y1} and
\cite{Y3},\,pp.\,65-70 for a brief detail. The above isomorphisms are
compatible with the action of the Hecke operators. Moreover,
according to the work of Skoruppa and Zagier \cite{SZ}, there is a
Hecke-equivariant correspondence between Jacobi forms of weight
$k$ and index $m$, and certain usual modular forms of weight
$2k-2$ on $\G_0(N).$

\psb \indent
Now we give the definition of Heegner
points of an elliptic curve $E$ over $\BQ$. By \cite{BCDT}, $E$ is
modular and hence one has a surjective holomorphic map
$\phi_E:X_0(N)\lrt E(\BC).$ Let $K$ be an imaginary quadratic
field of discriminant $D$ such that every prime divisor $p$ of $N$
is split in $K$. Then it is easy to see that $(D,N)=1$ and $D$ is
congruent to a square $r^2$ modulo $4N$. Let $\Theta$ be the set
of all $z\in\BH$ satisfying the following conditions
$$az^2+bz+c=0,\quad a,b,c\in\BZ,\ N|a,$$ $$b\equiv
r\,\,(\text{mod}\,2N),\qquad D=b^2-4ac.$$ Then $\Theta$ is
invariant under the action of $\G_0(N)$ and $\Theta$ has only a
$h_K$ $\G_0(N)$-orbits, where $h_K$ is the class number of $K$.
Let $z_1,\cdots,z_{h_K}$ be the representatives for these
$\G_0(N)$-orbits. Then $\phi_E(z_1),\cdots,\phi_E(z_{h_K})$ are
defined over the Hilbert class field $H(K)$ of $K$, i.e., the
maximal everywhere unramified extension of $K$. We define the
Heegner point $P_{D,r}$ of $E$ by
\begin{equation}
P_{D,r}=\sum_{i=1}^{h_K}\phi_E(z_i).
\end{equation} We observe that
$\epsilon=1$, then $P_{D,r}\in E(\BQ)$.

\psb \indent
Let $E^{(D)}$ be the
elliptic curve (twisted from $E$) given by
\begin{equation}
E^{(D)}\ :\ \
Dy^2=f(x).
\end{equation} Then one knows that the $L$-series of $E$ over
$K$ is equal to $L(E,s)\,L(E^{(D)},s)$ and that $L(E^{(D)},s)$ is
the twist of $L(E,s)$ by the quadratic character of $K/\BQ$.

\begin{theorem}[Gross-Zagier\,\cite{GZ, BCDT}]
	Let $E$ be an
	elliptic curve over $\BQ$ of conductor $N$ such that
    $\epsilon =-1.$ Assume
	that $D$ is odd. Then
	\begin{equation}
	L'(E,1)\,L(E^{(D)},1)=c_E\,u^{-2}\,|D|^{-{\frac 12}}\,{\hat
		h}(P_{D,r}),
	\end{equation}where $c_E$ is a positive constant not
	depending on $D$ and $r,\ u$ is a half of the number of units of
	$K$ and ${\hat h}$ denotes the canonical height of $E$.
\end{theorem}
\vskip 1mm
 Since
$E$ is modular by \cite{BCDT}, there is a cusp form of weight 2 with
respect to $\G_0(N)$ such that $L(f,s)=L(E,s).$ Let $\phi(z,w)$ be
the Jacobi form of weight 2 and index $N$ which corresponds to $f$
via the Skoruppa-Zagier correspondence. Then $\phi(z,w)$ has a
Fourier series of the form (6.2).
\psb \indent
B. Gross, W. Kohnen and D. Zagier\,\cite{GKZ} obtained the following result.
\begin{theorem}[Gross-Kohnen-Zagier\,\cite{GKZ, BCDT}]
Let $E$ be an elliptic curve over $\BQ$ with conductor $N$ and suppose that
$\epsilon=-1,\ r=1.$ Suppose that $(D_1,D_2)=1$ and $D_i\equiv
r_i^2\,(\text{mod}\,4N)\,(i=1,2).$ Then
$$L'(E,1)\,c((r_1^2-D_1)/(4N),r_1)\,c((r_2^2-D_2)/(4N),r_2)\,=\,c_E'<P_{D_1,r_1},
P_{D_2,r_2}>,$$ where $c_E'$ is a positive constant not depending
on $D_1,\,r_1$ and $D_2,\,r_2$ and $<\ ,\ >$ is the height pairing
induced from the N{\' e}ron-Tate height function ${\hat h}$, that
is,  ${\hat h}(P_{D,r})=<P_{D,r},P_{D,r}>$.
\end{theorem}
\vskip 1mm
We see from the
above theorem that the value of $<P_{D_1,r_1}, P_{D_2,r_2}>$ of
two distinct Heegner points is related to the product of the
Fourier coefficients
$c((r_1^2-D_1)/(4N),r_1)\,c((r_2^2-D_2)/(4N),r_2)$ of the Jacobi
forms of weight 2 and index $N$ corresponded to the eigenform $f$
of weight 2 associated to an elliptic curve $E$. We refer to \cite{Y4}
and \cite{Z} for more details.
\pmm\noindent
{\bf Corollary.} There is a point $P_0\in E(\BQ)\otimes_{\BZ}\BR$
such that $$P_{D,r}=c((r^2-D)/(4N),r)P_0$$ for all $D$ and
$r\,(D\equiv r^2\,(\text{mod}\,4N))$ with $(D,2N)=1.$ \psb The
corollary is obtained by combining Theorem H and Theorem I with
the Cauchy-Schwarz inequality in the case of equality.
\begin{remark}
R. Borcherds \,\cite{B} generalized the
Gross-Kohnen-Zagier theorem to certain more general quotients of
Hermitian symmetric spaces of high dimension, namely to quotients
of the space associated to an orthogonal group of signature $(2,b)$
by the unit group of a lattice. Indeed he relates the Heegner
divisors on the given quotient space to the Fourier coefficients
of vector-valued holomorphic modular forms of weight
$1+{\frac b2}$.
\end{remark}

\vskip 10mm
\section{Brief Reviews on the Works of Bhargava's School}
\vskip 0.25cm \ \ \
In this section, we briefly describe the recent works done by Bhargava's School. First we review the Selmer group.
Let $A=A(\overline{K})$ and $B=B(\overline{K})$ be abelian varieties over number field $K$ and let $f:A \rightarrow B$ be a nonzero isogeny with finite kernel
\begin{equation*}
A[f]=\{a\in A\;|\; f(a)=0\}.
\end{equation*}
Then we get a short exact sequence:
\begin{equation*}
\begin{tikzcd}
0 \ar[r]  & A[f]\ar[r, "\alpha"]& A \ar[r, "f"] & B\ar[r] & 0 .
\end{tikzcd}
\end{equation*}

Let $G_{K}=$Gal$(\overline{K}/K)$ be the Galois group of $\overline{K}$ over $K$. Then we have the following long exact sequence of Galois cohomology groups\;:
\begin{equation}
\begin{tikzcd}
0 \ar[r] & A[f]^{G_{K}}=A(K)[f] \ar[r] & A^{G_{K}}=A(K) \ar[r,"f"] & B^{G_{K}}=B(K) \ar[out=-30, in=150,"\delta_{1}"']{dll} \\
&H^1(G_{K},A[f])\rar   & H^1(G_{K},A) \rar   &\cdots  . &
\end{tikzcd}
\end{equation}
From (7.1) we obtain the following short exact sequence:
\begin{equation}
\begin{tikzcd}
0 \ar[r]  & B(K)/(f(A(K)))\ar[r, "\delta"]&H^1(G_{K},A[f]) \ar[r] &H^1(G_{K},A)[f]\ar[r] & 0.
\end{tikzcd}
\end{equation}
We let $K_v$ and $G_v$ the completion at $v$ and the decomposition group of $K_v$ respectively. We put $A_v:=A(\overline{K}_v)$.
Since $G_v$ acts on $A_v,\;B_v$ and $A_v[f]$, we get the short exact sequence:
\begin{equation}
\begin{tikzcd}
0 \ar[r]  & B_{v}(K_{v})/(f(A_{v}(K_{v})))\ar[r, "\delta"]&H^1(G_{v},A_{v}[f]) \ar[r] &H^1(G_{v},A_{v})[f]\ar[r] & 0.
\end{tikzcd}
\end{equation}

From the above short exact sequences (7.2) and (7.3), we have the following commutative diagram:
\begin{equation*}
\begin{tikzcd}[column sep=5mm]
0 \ar[r]  & B(K)/(f(A(K)))\ar[d]\ar{r}[name=U]{\delta} &H^1(G_{K},A[f])\ar[d,"\prod_{v}Res_{v}"] \ar{r}[name=X]{}&H^1(G_{K},A)[f]\ar[r]\ar[d,"\prod_{v}Res^{\ast}_{v}"] & 0 \\
0 \ar[r]  &\displaystyle\prod_{v} B_{v}(K_{v})/(f(A_{v}(K_{v})))\ar{r}[name=D]{\delta}&\displaystyle\prod_{v}H^1(G_{v},A_{v}[f]) \ar{r}[name=Y]{} &\displaystyle\prod_{v}H^1(G_{v},A_{v})[f]\ar[r] & 0
\ar[to path={(D) node[above,scale=2] {$\circlearrowleft$}  (U)}]{}\ar[to path={(Y) node[above,scale=2] {$\circlearrowleft$}  (X)}]{}.
\end{tikzcd}
\end{equation*}
Here $\mathrm{Res}_{v}(\xi)=\xi|_{G_{v}}$ for $\xi\in H^1(G_{K},A[f])$.
\begin{definition}
With the above notations, the \textsf{$f$-Selmer group of $A/K$} is defined by
\begin{eqnarray*}
\mathrm{Sel}^{(f)}(A/K)&=&\mathrm{Ker }\Big(H^1(G_{K},A[f])\rightarrow \prod_{v  }H^1(G_{v},A_{v})[f]\Big),\\
\hspace{-15mm}
&=&\bigcap_{v  }\mathrm{Ker }\Big(H^1(G_{K},A[f])\rightarrow H^1(G_{v},A_{v})[f]\Big).
\end{eqnarray*}
\end{definition}
\begin{definition}
The \textsf{Shafarevich-Tate group of $A/K$} is defined by
\begin{eqnarray*}
\Sha(A/K)&=&\mathrm{Ker }\Big(H^1(G_{K},A)\rightarrow \prod_{v  }H^1(G_{v},A_{v})\Big),\\
&=&\bigcap_{v  }\mathrm{Ker }\Big(H^1(G_{K},A)\rightarrow H^1(G_{v},A_{v})\Big).
\end{eqnarray*}
\end{definition}
\begin{theorem}
With the above notations, we get the following facts:
\begin{itemize}
\item[(a)] There is an exact sequence
\begin{equation*}
\begin{tikzcd}
0 \ar[r]  & B(K)/(f(A(K)))\ar[r]&\mathrm{Sel}^{(f)}(A/K) \ar[r] &\Sha(A/K)[f]\ar[r] & 0.
\end{tikzcd}
\end{equation*}
\normalsize
\item[(b)] The Selmer group $\mathrm{Sel}^{(f)}(A/K)$ is finite.
\end{itemize}
\end{theorem}

\begin{example}Let $A=B=E$ be an elliptic curve over $K=\mathbb{Q}$ and let $f=[m]$ be the multiplication by $m$ endomorphism.
Then we obtain the following exact sequence\;:
\begin{equation*}
\begin{tikzcd}
0 \ar[r]  & E(\mathbb{Q})/mE(\mathbb{Q})\ar[r]&\mathrm{Sel}^{[m]}(E/\mathbb{Q}) \ar[r] &\Sha(E/\mathbb{Q})[m]\ar[r] & 0
\end{tikzcd}.
\end{equation*}
$\mathrm{Sel}_{m}:=\mathrm{Sel}^{[m]}(E/\mathbb{Q})$ is called the $m$-Selmer group.
\end{example}
 Any elliptic curve $E$ over $\mathbb{Q}$ is isomorphic to the unique cubic curve $E_{A,B}$  in the plane of the form
\begin{equation*}
E_{A,B} : y^2 = x^3 + A x + B,
\end{equation*}
where $A,B \in \mathbb{Z}$ and for all primes $p: p^6\nmid A$  whenever $p^4\mid A$. Let $\mathscr{E}$ be the set of all such $E_{A,B}$. If $ E=E_{A,B}  \in \mathscr{E}$, then we define the (naive) height of $E$ by
\begin{equation*}
H(E_{A,B}) := \text{max}\lbrace 4 \vert A^3 \vert\;,\; 27 B^2 \rbrace.
\end{equation*}
For $X \in \mathbb{R}_{>0},$ we define $\mathscr{E}_{<X}:=\lbrace E\in \mathscr{E} \;\vert\; H(E)<X\rbrace.$
\begin{definition}
For any $\phi:\mathscr{E}\rightarrow\mathbb{R}$, we define
\begin{equation*}
\text{Average}_{\mathscr{E}}(\phi)\big(=\text{Aver}_{\mathscr{E}}(\phi)\big):=\lim_{X\to\infty}\dfrac{\sum_{ E\in \mathscr{E}_{<X}}\phi(E)}{\sum_{ E\in \mathscr{E}_{<X}}1}.
\end{equation*}
if the limit exists.
Define $\overline{\text{Aver}}_{\mathscr{E}}(\phi)$ and  $\underline{\text{Aver}}_{\mathscr{E}}(\phi)$ using limsup or liminf.
If the property $P$ can be identified with characteristic fucntion $\chi_P:\mathscr{E}\longrightarrow\lbrace 0,1 \rbrace$,
$$
\text{Prob}_{\mathscr{E}}(P):=\text{Average}(\chi_P).
$$
Similarly we define $\overline{\text{Prob}}_{\mathscr{E}}(P)$ and  $\underline{\text{Prob}}_{\mathscr{E}}(P)$.
\end{definition}
Bhargava and Shanker \cite{BS1,BS2,BS3,BS4} proved the follwing results:
\begin{theorem}[Bhargava and Shankar] \hfill
\begin{enumerate}
\item$\text{Aver}_{\mathscr{E}}(\#\mathrm{Sel}_2)=3.$ (cf.\,\cite{BS1})
\item$\text{Aver}_{\mathscr{E}}(\#\mathrm{Sel}_3)=4.$ (cf.\,\cite{BS2})
\item$\text{Aver}_{\mathscr{E}}(\#\mathrm{Sel}_4)=7.$ (cf.\,\cite{BS3})
\item$\text{Aver}_{\mathscr{E}}(\#\mathrm{Sel}_5)=6.$ (cf.\,\cite{BS4})
\end{enumerate}
\end{theorem}
\begin{corollary}[Bhargava and Shankar\,\cite{BS4}, 2013]
$$
\text{Aver}_{\mathscr{E}}(\text{rank})\leqq 0.885.
$$
\end{corollary}

\begin{theorem}[T.\,Dokchitser - V.\,Dokchitser\,\cite{DTDV}, 2010]	
	Let $E$ be an elliptic curve over $\mathbb{Q}$ and let $p$ be any prime.
	Let $s_p(E)$ and $t_p(E)$ denote the rank of the $p$-Selmer group of $E$ and the rank of $E(\mathbb{Q})[p]$, respectively.
	Then the quantity $r_p(E):=s_p(E)-t_p(E)$ is even if and only if the root
	number of $E$ is $+1$. Here we recall that $s_p(E)$ and $t_p(E)$
are defined as
\begin{equation*}
  s_p(E):= \log_p (\# \mathrm{Sel}_p (E)) \qquad {\rm and}\qquad
  t_p(E):= \log_p (\# E(\BQ)[p])
\end{equation*}
respectively.
\end{theorem}
Bhargava and Shanker \cite{BS4,BSk} proved the follwing:
\begin{theorem}\label{7thm08} \hfill
\begin{enumerate}
\item $\underline{\text{Prob}}_{\mathscr{E}}(r=0)\geqq 0.2062$.  (cf.\,\cite{BS4})
\item $\underline{\text{Prob}}_{\mathscr{E}}(r=1)\geqq 0.2612$.  (cf.\,\cite{BS4})
\item $\underline{\text{Prob}}_{\mathscr{E}}(r=0 \text{ or } r=1)\geqq 0.8375$. (cf.\,\cite{BS4})
\item $\underline{\text{Prob}}_{\mathscr{E}}(r=r_{\text{an}}=0)$ is positive. (cf.\,\cite{BS2})
\item $\underline{\text{Prob}}_{\mathscr{E}}(r=r_{\text{an}}=1)$ is positive. (cf.\,\cite{BSk})
\end{enumerate}
\end{theorem}
\noindent
The results (4) and (5) in Theorem \ref{7thm08} imply that a positive proportion of $\mathscr{E}$ satisfies the BSD conjecture.

\vskip 3.5mm
Quite recently Bhargava, Elkies and Shnidman \cite{BES} computed the average size of the $\phi_k$-Selmer group as $k$ varies over the integers, where
\begin{equation*}
  \phi_k:E_k \lrt E_{-27k},\qquad E_k:\ y^2=x^3+k
\end{equation*}
is a natural 3-isogeny.

\vskip 10mm
\section{Adjoint Orbits of $SL(2,\BR)$}
\vskip 0.25cm \ \ \
In this section, we describe the adjoint orbits of the special linear group $SL(2,\BR)$ explicitly.

\vskip 0.25cm
For brevity, we write $G=SL(2,\BR)$ and let $K=SO(2)$ be a maximal
compact subgroup of $G$. The Lie algebra $\fg$ of $G$ is given by

\begin{equation*}
\fg=\left\{ \begin{pmatrix} x & \ y \\ z & -x\end{pmatrix}\Big| \ x,y,z\in
\BR\,\right\}.
\end{equation*}
We put
\begin{equation*}
X=\begin{pmatrix} 1 & \ 0 \\ 0 & -1
\end{pmatrix},\quad
Y=\begin{pmatrix} 0 & 1 \\ 1 & 0
\end{pmatrix},\quad
Z=\begin{pmatrix} \ 0 & 1 \\ -1 & 0
\end{pmatrix}.
\end{equation*}
Then the set $\left\{ X,Y,Z\right\}$ forms a basis for $\fg$. We
define an element $F(x,y,z)\in \fg $ by
\begin{equation}
F(x,y,z):=xX+yY+zZ=\begin{pmatrix} x & y+z \\ y-z & -x
\end{pmatrix}.
\end{equation}
Then we have the relations
\begin{equation*}
X^2+Y^2-Z^2=3I_2,\quad [X,Y]=2Z,\quad
[X,Z]=2Y,\quad [Y,Z]=-2X.
\end{equation*}
\newcommand\al{\alpha}
\newcommand\Cal{\mathcal}
It is easy to see that $X$ and
$Y$ are hyperbolic elements and $Z$ is an elliptic element. For a
nonzero real number $\alpha$, the $G$-orbit of $\alpha X$ is represented
by the one-sheeted hyperboloid
\begin{equation}
x^2+y^2-z^2=\alpha^2. \label{3eq05}
\end{equation}
The $G$-orbit of $\al Y\,(\alpha\in \BR^{\times})$ is also represented by
the hyperboloid (\ref{3eq05}). The $G$-orbit of $\al Z\,(\al\in
\BR^{\times})$ is represented by two-sheeted hyperboloids
\begin{equation}
x^2+y^2-z^2=-\al^2.\label{3eq06}
\end{equation}
Since
$$F(x,y,z)^2=(x^2+y^2-z^2)\cdot I_2,$$ we have for any
$k\in\BZ^+,$ $$F(x,y,z)^{2k}=(x^2+y^2-z^2)^k\cdot I_2.$$
Thus we
see that $F(x,y,z)$ is nilpotent if and only if $x^2+y^2-z^2=0.$
Therefore the set ${\Cal N}_{\BR}$ of all nilpotent elements in
$\fg$ is given by

\begin{equation}
{\Cal N}_{\BR}=\left\{ F(x,y,z)=\begin{pmatrix} x &
y+z\\ y-z & -x\end{pmatrix} \Big|\ x^2+y^2-z^2=0\right\}.
\end{equation}

We put
\begin{equation}
S={\frac 12}(Y+Z)=\begin{pmatrix} 0 & 1\\ 0 & 0\end{pmatrix},\quad
T= {\frac 12}(Y-Z)=\begin{pmatrix} 0 & 0\\ 1 & 0\end{pmatrix}.
\end{equation}
Obviously $S$ and $T$ are nilpotent elements in ${\Cal N}_{\BR}$
and they satisfy
\begin{equation}
[X,S]=2S\,,\qquad[X,T]=-2T\,,\qquad[S,T]=X\,\label{3eq08}
\end{equation}
and
\begin{equation}
\theta(X)=-X\,,\qquad\theta(S)=-T,\qquad\theta(T)=-S. \label{3eq09}
\end{equation}
where $\theta$ is the Cartan involution defined by  $\theta(g)= -{}^tg$ for $g$ in $\mathfrak{g}$. \\
According to equation (\ref{3eq08}) and (\ref{3eq09}), $\{X,\,S,\,T\}$ and $\{-X,\,-S,\,-T\}$ are KS-triples in $\mathfrak{g}$.
\par The $G$-orbit of $\al S\,(\al\in \BR^{\times})$ is represented
by the cone
\begin{equation}
x^2+y^2-z^2=0,\quad (x,y,z)\neq (0,0,0)
\end{equation}
depending on the sign of $\al$.
\vskip 0.3cm
If $\al > 0,$ the $G$-orbit of $\al S$ is characterized by the one-sheeted
cone
\begin{equation}
x^2+y^2-z^2=0,\quad z>0.\label{3eq10}
\end{equation}

If $\al < 0,$ the
$G$-orbit of $\al S$ is characterized by the one-sheeted cone

\begin{equation}
x^2+y^2-z^2=0,\quad z<0.\label{3eq11}
\end{equation}
The $G$-orbits of $\al
T\,(\al>0)$ are characterized by the one-sheeted cone (\ref{3eq11}) and
the $G$-orbits of $\al T\,(\al<0)$ are characterized by the
one-sheeted cone (\ref{3eq10}).

\vskip 0.3cm
We define
the $G$-orbits ${\Cal N}_{\BR}^+$ and ${\Cal N}_{\BR}^-$ by
\begin{equation}
{\Cal N}_{\BR}^+=G\cdot S =
\left\{ gSg^{-1}\in {\mathcal N}_{\BR}\,|\ g\in G\,\right\}
\end{equation}
and
\begin{equation}
{\Cal N}_{\BR}^-=G\cdot T =
\left\{ gTg^{-1}\in {\mathcal N}_{\BR}\,|\ g\in G\,\right\}.
\end{equation}
Then we obtain

\begin{equation}
{\Cal N}_{\BR}={\Cal N}_{\BR}^+\cup \left\{ 0\right\}\cup {\Cal N}_{\BR}^-.\label{3eq13}
\end{equation}
According to (\ref{3eq05}),\,(\ref{3eq06}) and (\ref{3eq13}), we
see that there are infinitely many hyperbolic orbits and elliptic
orbits, and on the other hand there are only three nilpotent
orbits in $\fg$.

\vskip 10mm
\section{Final Remarks}
\vskip 0.25cm \ \ \
In the final section, we describe the relation between the nilpotent orbits of $SL(2,\BR)$ and CM points, and
propose several conjectures relating to the BSD conjecture.

\vskip 2mm
Let $G=SL(2,\BR)$ and let $\frak g$ be the Lie algebra of $G$. Let
$\mathcal E_\BR,\ \mathcal H_\BR$ and $\mathcal N_\BR$ be the set of
all $G$-elliptic orbits in $\frak g$, the set of
all $G$-hyperbolic orbits in $\frak g$, and the set of
all $G$-nilpotent orbits in $\frak g$ respectively.
Let $\exp:\frak g\longrightarrow G$ be the exponential map and let
$\phi:G\longrightarrow \BH$ be the map defined by
\begin{equation*}
  \phi (g)={{ai+b}\over {ci+d}},\qquad
  g=\begin{pmatrix}
      a & b \\
      c & d
    \end{pmatrix}\in G,\quad i=\sqrt{-1}.
\end{equation*}
We define the composition map $\omega:=\phi\circ \exp :\frak g\longrightarrow \BH.$

\begin{proposition}
\begin{equation*}
 \omega ({\mathcal N_\BR})=\mathbb H .
\end{equation*}
\end{proposition}
\vskip 1mm\noindent
{\it Proof.}
The proof can found in \cite{KY}.
\hfill $\Box$

\vskip 2mm
First we recall the thirteen CM points
\begin{eqnarray*}
  \al_1 &=& { {1+i\,\sqrt{3}} \over 2},\qquad \al_2=i,\qquad
  \al_3= { {1+i\,\sqrt{7}} \over 2},\qquad \al_4=i\,\sqrt{2},\\
  \al_5 &=& { {1+i\,\sqrt{11}} \over 2}, \qquad
  \al_6=i \,\sqrt{3},\qquad \al_7= 2\,i,\qquad
  \al_8={ {1+i\,\sqrt{19}} \over 2}, \\
  \al_9 &=& { {1+i\,\sqrt{27}} \over 2}, \qquad \al_{10}=i\,\sqrt{7},\qquad
  \al_{11}= { {1+i\,\sqrt{43}} \over 2},\\
  \al_{12} &=& { {1+i\,\sqrt{67}} \over 2},
        \qquad \al_{13}={\ {1+i\,\sqrt{163}} \over 2}.
\end{eqnarray*}
We put $q:=e^{2\pi i\tau}\ (\tau\in\BH).$
The $q$-expansion of the modular invariant $j(\tau)$
is given by
$$j(\tau)=q^{-1}+744+\sum_{n=1}^\infty c(n)q^n=q^{-1}+744+196884q+21493760q^2+
\cdots.$$
It is known that
\begin{eqnarray*}
  j(\al_1) &=&0,\qquad  j(\al_2)=2^6\!\cdot\! 3^3,
  \qquad  j(\al_3)=-3^3\!\cdot\!5^3,\qquad
  j(\al_4)=2^6\!\cdot\!5^3,\\
  j(\al_5) &=& -2^{15}, \qquad
  j(\al_6)=2^4\!\cdot\!3^3\!\cdot\!5^3,
  \qquad  j(\al_7)=2^3\!\cdot\!3^3\!\cdot\!11^3,
  \qquad  j(\al_8)=-2^{15}\!\cdot\!3^3,  \\
  j(\al_9) &=& -2^{15}\!\cdot\!3\!\cdot\!5^3, \qquad
  j(\al_{10})= 3^3\!\cdot\!5^3\!\cdot\!17^3,  \qquad
  j(\al_{11})=-2^{18}\!\cdot\! 3^3\!\cdot\!5^3,\\
  j(\al_{12}) &=& -2^{15}\!\cdot\! 3^3\!\cdot\!5^3\!\cdot\! 11^3,
\qquad
  j(\al_{13})=-2^{18}\!\cdot\!3^3\!\cdot\!5^3\!\cdot\!
23^3\!\cdot\!29^3.
\end{eqnarray*}

\vskip 0.25cm
For any element $\tau\in\BH$, we let
\begin{equation*}
  E_\tau:= \BC/\Lambda_\tau
\end{equation*}
be the elliptic curve over $\BC$,
where $\Lambda_\tau:=\BZ+\BZ\tau$ is the lattice in $\BC$.

\vskip 3mm
For any $g=\begin{pmatrix}
             a & b \\
             c & d
           \end{pmatrix}\in G$, we write
\begin{equation*}
  g<\!\tau\!>:={{a\tau +b}\over {c\tau +d}},\qquad \tau\in\BH.
\end{equation*}
Let
\begin{equation*}
  \mathcal O_j:=K\!<\!\alpha_j\!>=\left\{ k<\!\alpha_j\!>\,|\ k\in K\,\right\},\qquad 1\leq j\leq 13
\end{equation*}
be the $K$-orbit of $\alpha_j$. Clearly $\mathcal O_2=\{ i \}.$

\vskip 3mm
We propose the following conjecture.

\begin{conjecture}
Let $\tau$ be an element of $\mathcal O_j\ (1\leq j\leq 13)$
such that $E_\tau$ is defined over $\BQ$. Then the BSD conjecture for $E_\tau$ holds.
\end{conjecture}

For any $\alpha_j\ (j\neq 2,\ 1\leq j\leq 13)$, by Proposition 9.1,
we can choose $F(a_j,b_j,c_j)\in \mathcal N_\BR$ with $c_j\neq 0.$
We define
\begin{equation*}
  \mathfrak C_j:=\left\{ F(x,y,c_j)\in \mathcal N_\BR\,|\
  x^2+y^2=c_j^2\,\right\},\qquad j\neq 2,\ 1\leq j\leq 13
\end{equation*}
and
\begin{equation*}
  \mathfrak X_j:=\omega (\mathfrak C_j),\qquad j\neq 2,
  \ 1\leq j\leq 13.
\end{equation*}
We also define
\begin{equation*}
  \mathfrak D_j:=\left\{ F(x,b_j,z)\in \mathcal N_\BR\,|\
  x^2+b_j^2=z^2\,\right\},\qquad j\neq 2,\ 1\leq j\leq 13
\end{equation*}
and
\begin{equation*}
  \mathfrak Y_j:=\omega (\mathfrak D_j),\qquad j\neq 2,
  \ 1\leq j\leq 13.
\end{equation*}
Finally we define
\begin{equation*}
  \mathfrak E_j:=\left\{ F(a_j,y,z)\in \mathcal N_\BR\,|\
  a_j^2+y^2=z^2\,\right\},\qquad j\neq 2,\ 1\leq j\leq 13
\end{equation*}
and
\begin{equation*}
  \mathfrak Z_j:=\omega (\mathfrak E_j),\qquad j\neq 2,
  \ 1\leq j\leq 13.
\end{equation*}

\vskip 3mm
We propose the following conjectures.
\begin{conjecture}
Let $\tau$ be an element of
$\mathfrak X_j \,(j\neq 2,\,1\leq j\leq 13)$
such that $E_\tau$ is defined over $\BQ$. Then the BSD conjecture for $E_\tau$ holds.
\end{conjecture}

\begin{conjecture}
Let $\tau$ be an element of
$\mathfrak Y_j \,(j\neq 2,\,1\leq j\leq 13)$
such that $E_\tau$ is defined over $\BQ$. Then the BSD conjecture for $E_\tau$ holds.
\end{conjecture}

\begin{conjecture}
Let $\tau$ be an element of
$\mathfrak Z_j \,(j\neq 2,\,1\leq j\leq 13)$
such that $E_\tau$ is defined over $\BQ$. Then the BSD conjecture for $E_\tau$ holds.
\end{conjecture}


\begin{thebibliography}{99}
	
\bibitem{BS1} M. Bhargava and A. Shankar, {\em Binary quartic forms having bounded invariants, and the boundedness of the average rank of elliptic curves}, Annals of Mathematics
    {\bf 181}, no.\,1\,(2015), 191--242.	
	
\bibitem{BS2} M. Bhargava and Arul Shankar, {\em Ternary cubic forms having bounded invariants, and the existence of a positive proportion of elliptic curves having rank 0}, Annals of Mathematics {\bf 181}, no.\,2\,(2015), 587--621.

\bibitem{BS3} M. Bhargava and A. Shankar, {\em The average number of elements in the 4-Selmer groups of elliptic curves is 7},  arXiv:1312.7333v1 [math.NT] 27 Dec 2013.

\bibitem{BS4} M. Bhargava and A. Shankar, {\em The average size of the 5-Selmer group of elliptic curves is 6, and the average rank is less than 1}, arXiv:1312.7859v1 [math.NT] 30 Dec 2013.

\bibitem{BSk} M. Bhargava and C. Skinner, {\em A positive proportion of elliptic curves over have rank one}, arXiv:1401.0233v1 [math.NT] 1 Jan 2014.

\bibitem{BES} M. Bhargava, N. Elkies, and A. Shnidman, {\em The average size of the 3-isogeny Selmer groups of elliptic curves $ y^ 2= x^ 3+ k$}, arXiv:1610.05759v1 [math.NT] 18 Oct 2016.

\bibitem{BSD1} B. Birch and H.P.F. Swinnerton-Dyer, {\em Notes on elliptic curves  (I) }, J. Reine Angew. Math. {\bf 212} (1963), 7--25.

\bibitem{BSD2} B. Birch and H.P.F. Swinnerton-Dyer, {\em Notes on elliptic curves  (II) }, J. Reine Angew. Math. {\bf 218} (1965), 79--108.
\bibitem{B} R. Borcherds, {\em The Gross-Kohnen-Zagier
theorem in higher dimensions  }, Duke Math. J. {\bf 97}, no.\,2 (1999), 219--233.


\bibitem{BCDT} C. Breuil, B. Conrad, F. Diamond and R. Taylor, {\em On the modularity of elliptic curves over $\BQ$ }, Journal of AMS {\bf109} (2001), 843--939.

\bibitem{BFH} B. Bump, S. Friedberg and J. Hoffstein, {\em
Nonvanishing theorems for $L$-functions of modular forms and their
derivatives }, Invent. Math. {\bf 102} (1990), 543--618.

\bibitem{CZ1} H.-D. Cao and X.-P. Zhu, {\em Complete Proof of the Poincaré and Geometrization Conjectures-application of the Hamilton-Perelman theory of the Ricci flow}, Asian Journal of Mathematics {\bf10.2} (2006), 165-492.

\bibitem{CZ2} H.-D. Cao and X.-P. Zhu, {\em Erratum to "A complete proof of the Poincaré and geometrization conjectures — application of the Hamilton–Perelman theory of the Ricci flow", Asian J. Math., Vol. 10, No. 2, 165--492, 2006}, Asian Journal of Mathematics. {\bf 10(4)}, 663–-664.

\bibitem{CW} J. Coates and A. Wiles, {\em On the Birch-Swinnerton-Dyer conjecture} Invent. Math. {\bf 39} (1977), 223--252.

\bibitem{D} P. Deligne, {\em Formes modulaires et repr{\'e}sentations $\ell$-adiques}, S{\'e}minaire Bourbaki. Vol. 1968/69: Expos{\'e}s 347--363, Exp. No. 355, 139--172,
    Lecture Notes in Math., {\bf 175}, Springer, Berlin (1971).

\bibitem{DS} P. Deligne and J.-P. Serre, {\em Formes modulaires de poids 1}, Ann. Sci. Ec. Norm. Sup., IV. Ser.\,{\bf 7} (1974), 507--530.

\bibitem{DTDV} T. Dokchitser - V. Dokchitser, {\em On the Birch-Swinnerton-Dyer quotients modulo squares}, Annals of Mathematics {\bf 172}, no.\,1\,(2010), 567--596.

\bibitem{EZ} M. Eichler and D. Zagier, {\em The theory of Jacobi forms}, Birkh{\"a}user {\bf 55} (1985).

\bibitem{GZ} B. Gross and D. Zagier, {\em Heegner points
and derivatives of $L$-series}, Invent. Math. {\bf 84} (1986), 225--320.

\bibitem{GKZ} B. Gross, W. Kohnen and D. Zagier, {\em
Heegner points and derivatives of $L$-series. II}, Math. Ann.
,{\bf 278} (1987), 497-562.

\bibitem{KL} B. Kleiner and J. Lott, {\em Notes on Perelman’s papers}, Geom. Topol {\bf 12.5}(2008), 2587--2855.	


\bibitem{Koh} W. Kohnen, {\em Modular forms of half
integral weight on $\G_0 (4)$}, Math. Ann. {\bf 248} (1980), 249--266.

\bibitem{K1} V. A. Kolyvagin, {\em Finiteness of $E (\BQ)$
and $\text{III}(E,\BQ)$ for a subclass of Weil curves  (Russian)
}, Izv. Akad. Nauk SSSR Ser. Mat. {\bf 52} (1988),
522--540, 670--671\,; English translation in Math. USSR-IZv. {\bf 32}\ (1980),
523--541.

\bibitem{K2}  V. A. Kolyvagin, {\em Euler systems, the Grothendieck
Festschrift  (vol. II), edited by P. Cartier and et al },
Birkh{\"a}user {\bf 87} (1990), 435-483.


\bibitem{KY} Y.-J. Kwon and J.-H. Yang, {\em Nilpotent Orbits of $SL(2,\mathbb R)$ and Elliptic Curves with Complex Multiplication}, in preparation.

\bibitem{M1} H. Maass, {\em {\"U}ber eine Spezialschar von
Modulformen zweiten Grades I }, Invent. Math. {\bf 52} (1979), 95--104.

\bibitem{M2} H. Maass, {\em {\"U}ber eine Spezialschar von
Modulformen zweiten Grades II}, Invent. Math. {\bf 53} (1979), 249--253.

\bibitem{M3} H. Maass, {\em {\"U}ber eine Spezialschar von
Modulformen zweiten Grades III}, Invent. Math. {\bf 53} (1979), 255--265.


\bibitem{Ma1} B. Mazur, {\em Modular curves and the
Eisenstein series}, Publ. IHES {\bf 47} (1977), 33-186.

\bibitem{Ma2} B. Mazur, {\em Number Theory as Gadfly},
Amer. Math. Monthly {\bf 98} (1991), 593--610.

\bibitem{Mo} L. Mordell, {\em On the rational solutions of the indeterminate equations of the third and fourth degrees,} Proc. Cambridge Philos. Soc.(1922), pp. 179--192.

\bibitem{MT1} J. Morgan, G. Tian, {\em Completion of Perelman's proof of the geometrization conjecture}, arXiv:0809.4040

\bibitem{MT2} J. Morgan, G. Tian, {\em Ricci flow and the Poincaré conjecture}, Clay Mathematics Monographs {\bf 3}, Amer. Math. Soc. (2007).

\bibitem{MM} M.\,R. Murty and V.\,K. Murty, {\em Mean values
of derivatives of modular $L$-series }, Ann. Math. {\bf 133} (1991), 447--475.

\bibitem{P1} G. Perelman, {\em The entropy formula for the Ricci flow and its geometric applications}, arXiv:math/0211159v1
    [math.DG] 11 Nov 2002.

\bibitem{P2} G. Perelman, {\em Ricci flow with surgery on three-manifolds}, arXiv:math/0303109v1 [math.DG] 10 Mar 2003.

\bibitem{P3} G. Perelman, {\em Finite extinction time for the solutions to the Ricci flow on certain three-manifolds}, arXiv: math/0307245v1 [math.DG] 17 Jul 2003.

\bibitem{R} K. Rubin, {\em Elliptic curves with complex
multiplication and the BSD conjecture }, Invent. Math. {\bf 64}
 (1981), 455-470.

\bibitem{Sh} G. Shimura, {\em An $\ell$-adic method in the theory of automorphic forms}, 1968 (unpublished).

\bibitem{Si} J. Silvermann, {\em The Arithmetic of
Elliptic Curves},Springer-Verlag,  Graduate Text in Math.
{\bf 106} (1986).

\bibitem{SZ} N.-P. Skoruppa and D. Zagier, {\em Jacobi forms
and a certain space of modular forms }, Invent. Math. {\bf 94}
 (1988), 113--146.

\bibitem{Y1} J.-H. Yang, {\em Remarks on Jacobi forms of
higher degree }, Proceedings of the 1993 Conference on
Automorphic Forms and Related Topics, edited by J.-W. Son and
J.-H. Yang, Pyungsan Institute for Mathematical Sciences {\bf 1}
(1993), 33--58.

\bibitem{Y2} J.-H. Yang, {\em Note on Taniyama-Shimura-Weil
conjecture }, Proceedings of the 1994 Conference on Number
Theory and  Related Topics, edited by J.-W. Son and J.-H. Yang,
Pyungsan Institute for Mathematical Sciences {\bf 2} (1995
), 29--46.

\bibitem{Y3} J.-H. Yang, {\em Kac-Moody algebras, the Monstrous
Moonshine, Jacobi Forms and Infinite Products}, Proceedings of
the 1995 Symposium on Number Theory, Geometry and Related Topics,
edited by J.-W. Son and J.-H. Yang, Pyungsan Institute for
Mathematical Sciences {\bf 3} (1996), 13--82.

\bibitem{Y4} J.-H. Yang, {\em Past twenty years of the theory of
elliptic curves  (Korean)}, Comm. Korean Math. Soc.  {\bf 14}
 (1999), 449--477.

\bibitem{Y5} J.-H. Yang, {\em The Birch-Swinnerton-Dyer Conjecture}, Proceedings of the 2002 International Conference on Related Subjects to Clay Problems, the Institute of Pure and Applied Mathematics, Cheonbuk National Univ., Jeonju, Korea, Vol. {\bf 1} (2002), 134--151 or arXiv:math/0611423v1 [math.HO] 14 Nov 2006.

\bibitem{Z} D. Zagier, {\em $L$-series of Elliptic
Curves, the BSD Conjecture, and the Class Number Problem of Gauss
}, Notices of AMS {\bf 31} (1984), 739--743.

\end{thebibliography}
\end{document}